\documentclass[12pt]{article}

\usepackage{amsmath, amssymb, latexsym, amscd} 

\usepackage[dvips]{graphicx}
\usepackage{subcaption}

\newtheorem{theorem}{Theorem}[section]
\newtheorem{corollary}[theorem]{Corollary}

\newtheorem{statement}[theorem]{Statement}

\newtheorem{definition}[theorem]{Definition}
\newtheorem{remark}[theorem]{Remark}

\def\diam{\operatorname{diam}}

\def\aff{\operatorname{Aff}}
\def\dem{\operatorname{dem}}

\begin{document}


\title{A new simple family of Cantor sets in $\mathbb R^3$
all of whose projections are one-dimensional}
\author{Olga Frolkina\footnote{Supported by Russian Foundation of Basic Research; Grant
No.~19--01--00169.}\\
Chair of General Topology and Geometry\\
Faculty of Mechanics and Mathematics\\
M.V.~Lomonosov Moscow State University\\
and\\
Moscow Center for Fundamental and Applied Mathematics\\
Leninskie Gory 1, GSP-1\\
Moscow 119991, Russia\\
E-mail: olga-frolkina@yandex.ru
}

\maketitle

\begin{abstract}
In 1994, J.~Cobb
described a Cantor set in $\mathbb R^3$
each of whose projections into $2$-planes is
one-dimensional.
A series of works by other authors developing this field followed.
We present another very simple series of Cantor sets in $\mathbb R^3$
all of whose projections are connected and one-dimensional.
These are self-similar Cantor sets which go back to the work of Louis Antoine,
and we celebrate their centenary birthday in 2020--2021.


\end{abstract}

\section{Introduction}

In 1994, J.~Cobb
constructed a Cantor set in $\mathbb R^3$
such that its projection into every $2$-plane is
one-dimensional \cite{Cobb-projections}.
In the same paper, he stated a general question:
``Given $n>m>k>0$, does there exist a Cantor set in $\mathbb R^n$ such that
each of its projections
into $m$-planes is exactly $k$-dimensional?'' 
Such sets will be called $(n,m,k)$ Cantor sets, or briefly $(n,m,k)$-sets.

The case $m=k$ was already solved.
L.~Antoine
con\-struc\-ted 
a Cantor set in $\mathbb R^2$
whose projections coincide
with those of a regular hexagon
\cite[{\bf 9}, p.~272; and fig.~2 on p.~273]{Antoine-FM}.
(Another $(2,1,1)$ set is described in 
\cite[p.~124]{Cobb-projections};
for a different example and further results on planar sets with prescribed projections,
see \cite{Dijkstra-van-Mill}.) 
K.~Borsuk constructed an $(n,m,m)$ Cantor set in \cite{Borsuk};
as a corollary, he obtained a positive answer to the question of R.H.~Fox:
there exists a simple closed curve in 
$\mathbb R^3$ such that each of its plane projections is two-dimensional.

Later, examples for the following triples appeared:
for $(n,m,m-1)$ in \cite{Frolkina-proj} and
for $(n,n-1,k)$ in \cite{BDM} --- both generalized Cobb's idea;
and 
a different construction for $(n,n-1,n-2)$ in \cite{Frolkina-Greece-volume},
where further references can be found.
Interestingly, all projections of a \emph{typical} Cantor set are Cantor sets \cite{Frolkina-Arch-volume}.

Cobb's construction of a $(3,2,1)$ Cantor set is rather complicated;
the resulting set is tame. 
I suggested a different approach observing that
there exist Antoine's necklaces all of whose projections are one-dimensional
\cite[Thm. 1]{Frolkina-Greece-volume}.
To obtain a necklace with this property,
I  controlled the location and thickness of solid tori 
at prelimit levels of construction;
such control can also be generalized to all Cantor sets in all dimensions:
For each Cantor set $K\subset \mathbb R^n$, $n\geqslant 2$,
there is an arbitrarily small isotopy $\{ h_t \}: \mathbb R^n\cong \mathbb R^n$ 
such that
the projection of $h_1(K)$ into each  $(n-1)$-plane has dimension
$(n-2)$ \cite[Thm.~2]{Frolkina-Greece-volume}.

The control process is not trivial. 
It is based on the 
plank theorem 
(for the planar case, see
English translation of
A.~Tarski's and H.~Moese's papers (1931--32)
in \cite[Chapter 7]{MMS} where the history
of the problem is also discussed;
the high-dimensional case is treated in \cite{Bang}).

The control process is necessary. In fact:
For each Cantor set $K\subset \mathbb R^n$, $n\geqslant 2$,
there exists an arbitrarily small isotopy $\{ h_t \} : \mathbb R^n\cong \mathbb R^n$ 
such that 
for each proper subspace $\Pi \subset \mathbb R^n$
the dimension of the projection of $h_1(K)$ into $\Pi $ equals
$\dim \Pi $
\cite[Stat. 12]{Frolkina-Greece-volume}.

In this article, I restore the original definition of Antoine's necklaces 
from~\cite{Antoine-diss}. These sets are self-similar. 
In Theorem \ref{main}, we formulate
simple additional conditions 
which ensure that all projections of these sets are one-dimensional. Theorems \ref{chains} and \ref{suff} 
imply that these conditions can actually be satisfied.
Self-similarity eliminates the need for thickness control: 
one-dimensionality is now
derived from properties of
planar Lebesgue measure or Hausdorff measure. 

We also consider a high-dimensional case (Theorem \ref{measure}).
Using Rushing's idea to combine results of Shtan'ko and 
V\"{a}is\"{a}l\"{a} \cite[p.~599]{Rushing-wild},
we get a 
new shorter --- but less elementary ---  proof of the above proposition
concerning $(n,n-1,n-2)$-sets, 
see Corollary \ref{move}; 
the usage of the plank theorem is eliminated from the proof.
Finally, we observe that the sets 
$C^{s}$ from \cite{Rushing-wild} are $(3,2,1)$-sets assuming that
$1\leqslant s <2$ (Remark \ref{R}). 

Let us underline that the purpose of this note is to propose
a possibly simplest example of a $(3,2,1)$ Cantor 
set in a self-contained presentation. 
By this reason we first give an ``elementary'' exposition, and 
thereafter state some generalizations.

Finishing the introduction, I would like to recall some other 
results of L.~Antoine.

Antoine's necklace is a first example of a wild Cantor set in $\mathbb R^3$;
its description 
\cite[{\bf 78}, p.~91--92]{Antoine-diss},
briefly presented in \cite{Antoine},
is clear,
in contrast to vague descriptions
with chains of ``curved solid tori'' that are often found in later retelling.
(Note that in order to get a Cantor set, the diameters of prelimit tori should
tend to zero. For self-similar constructions, this holds automatically.
Constructions with e.g. ``four curved tori''  must be provided with a possibility proof.)
Antoine constructed
first wild arcs in $\mathbb R^3$.
His first example  \cite{Antoine},
\cite[{\bf 54--58}, p. 65--70]{Antoine-diss}
consists 
of an arc with a sequence of trefoils tied in it and convergent to the end point, 
united with a straight line segment.
The second wild arc \cite{Antoine}, \cite[{\bf 83}, p. 97]{Antoine-diss} contains an Antoine's necklace.
Antoine constructed an everywhere wild arc in $\mathbb R^3$ \cite{Antoine-FM}.
In \cite{Antoine-173},
Antoine announced
and in \cite{Antoine-FM}
described in detail 
an embedded $2$-sphere in $\mathbb R^3$
which contains 
an Antoine's necklace,
this is the first example of a wild surface; see \cite[Thm.~18.7]{Moise}.
I would 
recommend the whole volume
\cite{Antoine-100} for acquaintance with Antoine's destiny and results;
and especially
\cite{Ibisch-A}, \cite{Ibisch}, and \cite[Preface, Sec.~18, and p.~256]{Moise}
concerning the importance of his topological discoveries.

\subsection*{Notation and Conventions}

A plane in  $\mathbb R^3$ is a 2-dimensional affine subspace.

For a plane $\Pi \subset \mathbb R^3$,
the orthogonal projection map
$\mathbb R^3\to \Pi $
is denoted by $p_{\Pi }$.

For a set $X\subset \mathbb R^3$, its affine hull
is denoted by $\aff X$, and its diameter by $\diam X$.

For topological dimension theory, see \cite{HW};
for definition and properties of Lebesgue measure, refer to e.g. \cite{Bogachev}.
Definitions and basic properties 
of Hausdorff measure and Hausdorff dimension
can be found in \cite[Chap.~VII]{HW} of \cite[3.6]{DV}.
For Shtan'ko demension theory, see 
\cite{Shtanko-theory}; 
its presentation is given in \cite{Edwards} and \cite[3.4]{DV}.

\section{Main construction and results: ``elementary'' approach}

There are several (more or less restrictive)
definitions of Antoine's necklaces,
compare e.g.
\cite{Sher68}, \cite{GRZ2005}, \cite{Z}, \cite{Frolkina-Greece-volume}.
They depend on allowed shapes, numbers and positions of solid tori on prelimit stages.
For our purpose, we will only restrict ourselves by the class of self-similar 
Antoine's necklaces.
This construction comes back to L.~Antoine \cite{Antoine}, \cite{Antoine-diss}.

\begin{definition}\label{torus-r-R}
Let $\Pi $ be a plane in $\mathbb R^3$.
Let $D\subset \Pi $ be a disk of radius $r>0$
with center $Q$, and let $\ell \subset \Pi $ be a straight line such that $d(Q,\ell ) = R > r$.
A \emph{standard solid torus} $T$
is the solid torus of revolution 
 generated by
revolving $D$ in $\mathbb R^3$ about $\ell $.
The \emph{central circle} of $T$
is the circle generated by rotating the point $Q$. The \emph{center} of $T$ is the center
of its central circle.
For a given standard solid torus $T$, the corresponding radii
will be denoted by $r_T$ and $R_T$, the central circle by $C_T$, and the center by $Q_T$.
\end{definition}

\begin{figure}[!tbp]
  \begin{subfigure}[b]{0.5\textwidth}
  \centering
  \scalebox{0.6}{\includegraphics{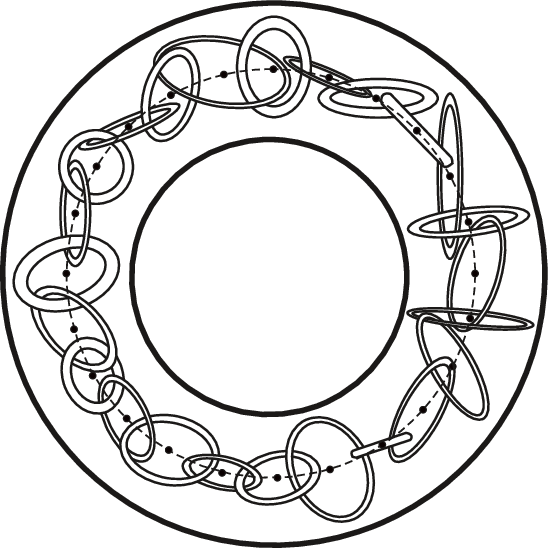}}
    \caption{A simple chain.}
    \label{fig:fig1a}
  \end{subfigure}
  \begin{subfigure}[b]{0.5\textwidth}
   \centering
  \scalebox{0.6}{\includegraphics{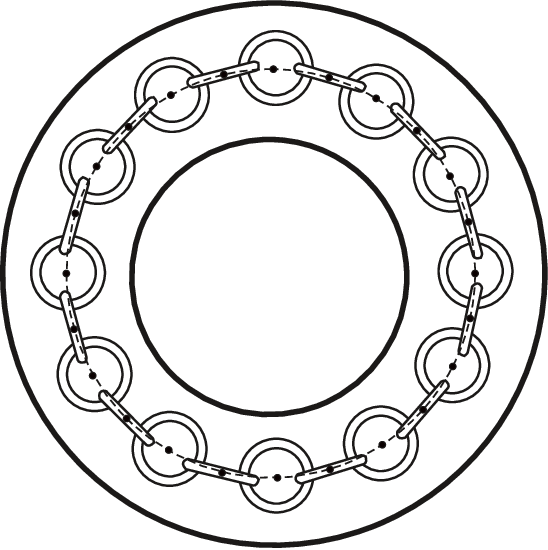}}
    \caption{A regular simple chain.}
    \label{fig:fig1b}
  \end{subfigure}
  \caption{Simple chains with $2m=24$ links.}
  \end{figure}

\begin{definition}\label{A-def}
\emph{A simple chain}
in a standard solid torus $T\subset \mathbb R^3$
is a finite family 
$T_1,\ldots , T_k$, $k\geqslant 3$, 
of pairwise disjoint 
standard solid tori such that

1) $T_1\cup\ldots\cup T_k\subset T$;

2) centers of $T_1,\ldots, T_k$
are subsequent vertices of a regular convex $k$-gon inscribed in the
central circle of $T$;

3) $T_i$ and $T_j$ are linked for
$|i - j |\equiv 1\mod k$, and are not linked
otherwise.

See Fig. \ref{fig:fig1a}.
\end{definition}

A simple chain $T_1,\ldots , T_k$ 
looks like a ``usual necklace'' which
winds once around the hole of $T$;
no one of the tori $T_i$ embraces the hole of $T$.

\begin{definition}\label{nice-chain}
A simple chain $T_1,\ldots , T_{2m}$, $m\geqslant 2$
in a standard solid torus $T\subset \mathbb R^3$
will be called 
\emph{regular}
if 

1) all $T_i$'s are congruent,

2)
$\aff C_{T_1} = \aff C_{T_3} = \ldots = \aff C_{T_{2m-1}} = \aff C_T$,

3)
$\aff C_{T_{2i}} \bot \aff C_{T}$ for $i=1,\ldots , m$, and

4)
each straight line
$\aff C_{T_{2i}} \cap \aff C_T$, $i=1,\ldots , m$,
touches the circle $C_T$.

See Fig. \ref{fig:fig1b}.
\end{definition}

\begin{figure}[ht]
\centering
\scalebox{1.4}
{\includegraphics{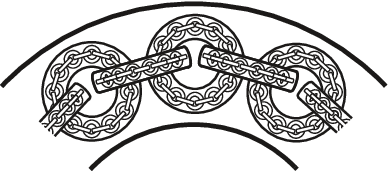}}
\caption{A piece of a self-similar Antoine's necklace: first stages of construction.}
\label{fig:fig2}
\end{figure}

\begin{definition}\label{self-sim}
Let $T\subset\mathbb  R^3$ 
be a standard solid torus, and
let $S_1, \ldots , S_k : \mathbb R^3\cong\mathbb R^3$
be similarity transformations such that 
$T_1 := S_1 (T) , \ldots , T_k:=S_k(T)$
is a simple chain in $T$.
For each $\lambda \in \mathbb N$ take
$$
M_{\lambda } =
\bigcup\limits_{(i_1,\ldots , i_\lambda ) \in \{ 1,\ldots , k\}^\lambda }
 S_{i_1} \circ S_{i_2}\circ \ldots \circ S_{i_\lambda } (T) .
 $$
 We have
$
T\supset M_1 \supset M_2 \supset M_3\supset \ldots 
$. The intersection
$\mathcal A  (T; S_1,\ldots , S_k) := \bigcap\limits_{\lambda =1}^\infty M_\lambda $
can be easily seen to be a Cantor set;
we call it \emph{a self-similar Antoine's Necklace generated by  
$(T; S_1,\ldots , S_k)$}.
For each $i=1,\ldots , k$ the piece
$  \mathcal A  (T; S_1,\ldots , S_k)  \cap T_i $ equals 
$S_i (\mathcal A  (T; S_1,\ldots , S_k))$, hence it
is geometrically similar to $\mathcal A  (T; S_1,\ldots , S_k)$.

A self-similar Antoine's Necklace will be called \emph{regular}
if the simple chain $T_1,\ldots , T_k$ in $T$ is regular.
(In particular, $k$ is even, and 
the similarity coefficients of all $S_i$'s are equal.)
See Fig. \ref{fig:fig2}.
\end{definition}

\begin{remark}
The necklace $\mathcal A (T; S_1, \ldots , S_{k}) $ 
is the invariant set
with respect to the set of similarities 
$\{ S_1, \ldots , S_{k} \}$;
we thus may omit $T$ from the notation.
For details, see \cite{Hutchinson}.
\end{remark}

L.~Antoine writes:
``Num\'erotons de~$1$ \`a~$k$ les $k$ tores d'ordre $\lambda + 1$ int\'erieurs
\`a un tore d'ordre $\lambda $.
Je supposerai que deux tores ayant des numeros cons\'ecutifs, ainsi que
les tores $k$ et~$1$, sont entlaces, alors que deux tores 
non cons\'ecutifs ne sont pas enlac\'es.
Je suppose enfin que les centres de ces $k$ tores sont sur la circonf\'erence d'ordre~$\lambda $ et qu'en joignant ces points dans l'ordre de leurs num\'eros,
on forme un polygone r\'egulier non \'etoil\'e.
Ces conditions sont possibles \`a r\'ealiser. Il suffit
pour cela que le rayon du m\'eridien d'un tore $T$ soit assez petit,
par rapport au rayon de la circonf\'erence $C$ de ce tore,
et que $k$ soit assez grand. On peut alors prendre pour tores d'ordre $\lambda +1$
des tores semblables au tore d'ordre $\lambda $. La possibilit\'e
de la construction a alors lieu quel que soit $\lambda $, $k$ \'etant fixe.
Pour simplifier quelques raisonnements ult\'erieurs, je
{\it supposerai que ce nombre $k$ est pair}''
\cite[{\bf 78}, p.~92]{Antoine-diss}.
Antoine's description does not explicitly indicate similarity transformations, in contrast to what is done in Definition \ref{self-sim}. But I believe that, taking similar tori, he had in mind precisely this. In \cite[I]{Antoine} and \cite{Antoine-FM} he carries out several more constructions using the idea of similarity.
Antoine assumed neither that chains have a special position in space,
nor he demanded that the solid tori which constitute a chain are congruent to each other.
To obtain examples with a greater degree of symmetry,
we added some restrictive conditions (regularity),
as is done in \cite{MR}, \cite{GRZ2005}, \cite{Z}, \cite{RS}.

Self-similar Antoine's necklaces with different $k$'s are inequivalently embedded,
this follows from
\cite[Thm.~2]{Sher68}.
For a regular self-similar Antoine's necklace, the inequality
$k \geqslant 20$ holds \cite{Z}.
A schematic picture of two inequivalently embedded 
self-similar Antoine's necklaces one of which is regular and the other is not is given 
in \cite{GRZ2005}; in that example, $k=60$.

\medskip

\begin{theorem}\label{main}
Let $\mathcal A (T; S_1, \ldots , S_{k}) \subset \mathbb R^3$ be a 
self-similar Antoine's necklace,
and let $s_i$ be the similarity coefficient of $S_i$ for $i=1,\ldots , k$.
If \ $s_1^2  + \ldots + s_k^2 <1$, 
then 
$p_{\Pi } \bigl( \mathcal A (T; S_1, \ldots , S_{k}) \bigr) $ is a connected one-dimensional set
for each plane $\Pi \subset \mathbb R^3$.

In particular, a regular self-similar Antoine's necklace 
 $\mathcal A (T; S_1, \ldots , S_{2m}) $
 with $2ms^2 < 1$ is a $(3,2,1)$-set;
 here $s$ is the similitude coefficient of $S_i$'s.
\end{theorem}

Below we will return to this result
from a more general point of view, see Theorem \ref{measure}.

\smallskip

Two points should be discussed.
The first one: do regular self-similar Antoine's necklaces exist?
(Equivalently, can one construct
a standard solid torus $T\subset \mathbb R^3$
together with a regular simple chain $T_1,\ldots , T_{2m}$ in it such that each
$T_i$ is similar to $T$~?)
If ``yes'', the second one: do regular self-similar Antoine's necklaces with $2ms^2 < 1$ exist?

If $T_1,\ldots , T_{2m}$ is a regular simple chain in $T$, 
and $T_i$'s are similar to $T$, then
the numbers $\frac{r_T}{R_T}$, $m$, and $s$ should satisfy several inequalities,
see \cite{Z}, \cite{RS}; it seems that no exact description of solutions exist.
A picture representing a regular self-similar Antoine's necklace with $2m=60$ can be found in \cite[p. 294]{GRZ2005}.
Using numerical methods,
it is shown in \cite{Z} that $2m\geqslant 20$ for any regular self-similar Antoine's necklace; and a picture for $2m = 20$ is given.
For $2m = 20,22,24,26,28$, existence of
regular self-similar Antoine's necklaces follows from
\cite[p.~146]{RS}, where computer software was also used.
Given a standard solid torus $T$
with $\frac{r_T}{R_T}$ sufficiently small,
it
seems plausible that for all sufficiently large integers $m$
there exists in $T$ 
a regular simple chain of $2m$ tori 
similar to $T$.
I have not find an exact statement in the literature.
More or less close assertions can be found without proofs 
in \cite[p. 274]{MR} (``It can be verified
by elementary calculations  that for all sufficiently big numbers $n$,  there exist
values of  $\lambda $ - not  too  big  and  not  too  small - and  corresponding thicknesses of $X$ such  that...''; formulaes and a picture follow saying that the standard 
solid torus $X$ contains a regular simple chain of $n$ tori similar to $X$ with the similarity coefficient $\lambda $)
and \cite[p. 132]{RS} (``...It seems clear that once this minimum is achieved, one can simply add more tori, in specific quantities derived below, to obtain another viable construction of the Antoine necklace'').
We give a rigorous statement and complete it by a purely mathematical 
proof.

\begin{theorem}\label{chains}
For a standard solid torus  $B  \subset\mathbb R^3$,
the following are equivalent:

(1) $R_B > 3r_B$;

(2) for each sufficiently large integer $m$
there exist a standard solid torus $T$
and a regular simple chain $T_1,\ldots , T_{2m}$ in $T$
such that each $T_i$ ($i=1,\ldots , 2m$) is congruent to $B$;

(3) for each sufficiently large integer $m$
there exist a standard solid torus $\mathbf T$
and a regular simple chain $T_1,\ldots , T_{2m}$ in $\mathbf T$
such that each $T_i$ ($i=1,\ldots , 2m$) is congruent to $B$,
and $\mathbf T$ is similar to $B$.
\end{theorem}

Now let us discuss 
the condition $2ms^2<1$ from Theorem \ref{main}.
Are we possible to satisfy it? Yes.
Briefly speaking, the similarity coefficient $s$ tends to zero as $m\to \infty $;
but the product $ms$ is bounded. Hence necklaces which satisfy the assumption
$2ms^2 < 1$ should exist.
We get several simple estimates:

\begin{theorem}\label{suff}
Let $\mathcal A (T; S_1, \ldots , S_{2m}) \subset \mathbb R^3$ be a 
regular self-similar Antoine's necklace,
and let $s$ be the similarity coefficient of $S_i$'s.

1) If \ 
$s < \frac{1}{2\pi }$, then
$2ms^2 <1$.

2)  If \ $\frac{r_T}{R_T} <  \frac{1}{2\pi - 1}  $, then $s < \frac{1}{2\pi }$.

3) If \ $2m\geqslant 40$, then $s < \frac{1}{2\pi }$.
\end{theorem}

\section{A generalization using Hausdorff measure}

\begin{definition}\label{tame-C}
A zero-dimensional compact set $X\subset \mathbb R^n$ is called \emph{tame}
if there exists
a homeomorphism $h$ of $\mathbb R^n$ 
onto itself such that 
$h(X)$ is a subset
of a straight line in $\mathbb R^n$;
otherwise, $X$ is called \emph{wild}.
\end{definition}

In $\mathbb R^2$ each zero-dimensional compactum
is tame 
\cite[{\bf 75}, p.~87--89]{Antoine-diss}
(one may also refer to
\cite[Chap.~13]{Moise}).

Antoine proved that the necklaces he described are wild
\cite[{\bf 80}, p.~93]{Antoine-diss}.
A later proof using fundamental group can
be found in \cite[Sect. 18]{Moise} or \cite[2.1]{DV}.

Antoine's construction was generalized to higher dimensions by W.A.~Blan\-kin\-ship and
by A.A.~Ivanov, see \cite{Frolkina-Greece-volume} for references;
hence 
$\mathbb R^n$ contains wild Cantor sets provided that $n\geqslant 3$.

\begin{theorem}\label{measure}
Let $X\subset \mathbb R^n$
be a wild Cantor set of a
zero $(n-1)$-dimensional Hausdorff measure:
$m_{n-1} (X) =0$.
(In particular, this holds if the Hausdorff dimension $\dim _H  X< n-1$.)
Then $X$ is an $(n,n-1,n-2)$-set.
\end{theorem}

For examples to this theorem, see Remark~\ref{R} and Corollary~\ref{move}.

\medskip

Theorem \ref{measure} covers Theorem \ref{main} except connectedness of projections
because of the following:

\begin{statement}\label{red}
Let $\mathcal A = \mathcal A (T; S_1, \ldots , S_{k}) \subset \mathbb R^3$ be a 
self-similar Antoine's necklace,
and let $s_i$ be the similarity coefficient of $S_i$ for $i=1,\ldots , k$.
If $s_1^2  + \ldots + s_k^2 <1$, then $m_2 (\mathcal A ) =0$.
\end{statement}

Statement \ref{red} follows from Moran's theorem,
see
\cite[p. 600]{Rushing-wild}, \cite[5.3(1)]{Hutchinson};
we give a short proof for our simple case. 

{\bf Proof of Statement 1.}
For each $\lambda \in \mathbb N$, we have $\mathcal A \subset M_\lambda $,
see Definition~\ref{self-sim}.
We can assume (zooming if necessary) that $\diam T = 1$.
Now $\mathcal A$ is covered by 
the union of $k^\lambda $ sets
$S_{i_1} \circ S_{i_2}\circ \ldots \circ S_{i_\lambda } (T) $, where
$ (i_1,\ldots , i_\lambda ) \in \{ 1,\ldots , k\}^\lambda $,
and
$$\diam  S_{i_1} \circ S_{i_2}\circ \ldots \circ S_{i_\lambda } (T) 
= s_{i_1} \cdot s_{i_2}\cdot \ldots \cdot s_{i_\lambda }  . $$
We have
$$ 
\sum\limits_{(i_1,\ldots , i_\lambda ) \in \{ 1,\ldots , k\}^\lambda }
\Bigl( \diam  S_{i_1} \circ S_{i_2}\circ \ldots \circ S_{i_\lambda } (T) 
 \Bigr)^2 \leqslant
$$
$$
\leqslant
\sum\limits_{(i_1,\ldots , i_\lambda ) \in \{ 1,\ldots , k\}^\lambda }
\bigl( s_{i_1} \cdot s_{i_2}\cdot \ldots \cdot s_{i_\lambda } \bigr) ^2 =
\bigl( s_1^2 +\ldots + s_k^2 \bigr) ^\lambda .
$$
 This tends to zero as $\lambda \to \infty $,
 hence the $2$-dimensional Hausdorff measure  $m_2  (\mathcal A) =0$.
 \hfill $\square $

\begin{remark}\label{R}
For each $s$ with $1\leqslant s \leqslant 3$, T.B.~Rushing \cite[Thm. 1, 2]{Rushing-wild}
described a wild Cantor set  $C^{s}\subset \mathbb R^3$
such that $\dim _H C^s = s$.
Theorem \ref{measure} implies that the sets $C^s$
are $(3,2,1)$ Cantor sets, provided that
$1\leqslant s <2$.
The sets $C^{s}$ are Antoine-type necklaces; 
their prelimit chains 
constructed under very careful control have a rather complex arrangement in space;
they do not satisfy our Definition \ref{self-sim}.
Rushing's scheme is too
complicated for our purpose;
the aim of this paper
is to give an example which is as simple as possible.
\end{remark}

\begin{remark}
It is apparently an open question, whether wild \emph{self-similar} Cantor sets in $\mathbb R^n$ for $n\geqslant 4$ do exist
\cite[p. 298, Question]{GRZ2005}.
T.B.~Rushing, without detailed proof, 
indicates the possibility of this \cite[p. 611]{Rushing-wild}.
\end{remark}

As another application of Theorem \ref{measure},
we get a new proof of our earlier result:

\begin{corollary}\label{move}\cite[Thm. 2]{Frolkina-Greece-volume}
Let $X\subset \mathbb R^n$, $n\geqslant 2$,
be any Cantor set.
For each $\varepsilon >0$
there exists an $\varepsilon $-isotopy
$\{ h_t \} :\mathbb R^n\cong \mathbb R^n$
such that $h_1(X)$ is an $(n,n-1,n-2)$-set.
\end{corollary}

The new proof is shorter but less elementary
since it refers to strong results of Shtan'ko, Edwards and
V\"{a}is\"{a}l\"{a}. 

{\bf Proof of Corollary 1.}
The case of a tame $X$ reduces to known results, 
see \cite{Frolkina-Greece-volume} for details.

Let $X$ be wild.
We start essentially as in Rushing's remark
\cite[p. 599]{Rushing-wild}.
Suppose we are given an arbitrary $\varepsilon >0$.
According to V\"{a}is\"{a}l\"{a}'s theorem, there exists an isotopy
$\{ h_t \} :\mathbb R^n\cong \mathbb R^n$
such that the Hausdorff dimension
$\dim _H (h_1(X))$ equals the embedding dimension $ \dem X$
 \cite{Vaisala}, \cite[Thm. 3.6.2]{DV}.
By the Edwards construction \cite[p.~208--209]{Edwards}, \cite[p.~168]{Vaisala}, 
this isotopy can be chosen to be an $\varepsilon $-isotopy.
For a wild Cantor set $X\subset \mathbb R^n$, we have $\dem X \leqslant n-2$ 
(for $n\geqslant 4$, see \cite[Remark~2, Prop.~1]{Shtanko-theory};
for $n=3$ use \cite[Lemma~4]{Armentrout}).
Consequently $\dim _H (h_1 (X)) \leqslant n-2$, so
 $h_1(X)$ is an $(n,n-1,n-2)$-set  by Theorem \ref{measure}.

(Actually, we have $\dem X = n-2$ for a wild Cantor set $X\subset\mathbb R^n$.
This is a special case of strong results from geometric topology due to several authors, see 
\cite{Shtanko-theory},
\cite[p.~195; Thm.~1.4 and its proof]{Edwards}, 
\cite[p.~598]{Rushing-wild}, 
\cite[Thm.~3.4.11, Exer.~3.4.1]{DV}, \cite[p.~5]{BDVW} for statements, 
historical details and references.)
\hfill $\square $

\section{Proofs of Theorems \ref{main}--\ref{measure}}

{\bf Proof of Theorem 1.}
We already derived Theorem \ref{main} (except connectedness
of projections) from Theorem
\ref{measure}.
Here we give an independent ``elementary'' proof.

For brevity, denote
$\mathcal A  = \mathcal A (T; S_1, \ldots , S_{k})$.
By
\cite[{\bf 80-82}, p. 93-96]{Antoine-diss},
no topological $2$-sphere separates a necklace.
Hence all projections of a necklace are connected (see
\cite[Stat. 2]{Frolkina-Greece-volume} for details).

Let $\Pi \subset \mathbb R^3$ be any plane.
It is clear that $p_{\Pi }  \mathcal A   $ contains more than one point;
so it remains to show that
$\dim p_{\Pi } \mathcal A \neq 2$.

A subset $X\subset  \mathbb R^n$ 
is $n$-dimensional iff it contains an $n$-ball
\cite[Thm. IV~3]{HW}.
Therefore it suffices to prove that
$\lambda _2 ( p_{\Pi } \mathcal A ) = 0$,
where $\lambda _2 $ denotes planar Lebesgue measure.
Suppose that $\lambda _2 (p_\Pi \mathcal A ) > 0$.
(In $\mathbb R^n$, $n$-dimensional Lebesgue and Hausdorff measures differ
by a constant multiple. We thus may use $m_2$ instead of $\lambda _2$.
Here I decided to use $\lambda _2$ due to its wider popularity.)

For each non-empty compact subset $X\subset \mathbb R^3$ 
define a real number
$$
\alpha (X) =
\sup \{
\lambda _2 (p_L X ) \ | \ 
L \text{ is a plane in } \mathbb R^3
\}.
$$
Any similarity transformation of $\mathbb R^n$
is a composition of
a scaling and an isometry.
As a corollary, for any
similarity transformation 
$S:\mathbb R^3\to \mathbb R^3$ and any non-empty compactum $X\subset \mathbb R^3$
we have
$\alpha (S(X)) = s^2 \cdot \alpha (X)$, where $s$ is the similarity coefficient of $S$.

By assumption, $\alpha (\mathcal A) >0$.
Take a positive number 
$$\varepsilon 
\leqslant  \alpha (\mathcal A) - \bigl( s_1^2 + \ldots + s_k^2 \bigr) \cdot \alpha  (\mathcal A) . $$
Let $ L(\varepsilon ) $ be any plane such that
 $
\alpha (\mathcal A ) - \varepsilon < \lambda _2(p_{L(\varepsilon )} \mathcal A )
\leqslant \alpha (\mathcal A) $.
Recall that
$$
 \mathcal A = S_1 ( \mathcal A ) \cup \ldots \cup S_{k}( \mathcal A ) . $$
We have
$$
\alpha (\mathcal A) - \varepsilon < \lambda _2(p_{L(\varepsilon )} \mathcal A ) =
\lambda _2\Bigl( p_{L (\varepsilon )} \bigl( S_1 (\mathcal A ) \cup
\ldots \cup S_{k}( \mathcal A ) \bigr) \Bigr) \leqslant
$$
$$
\leqslant
\lambda _2\Bigl( p_{L (\varepsilon )} \bigl( S_1 (\mathcal A ) \bigr) \Bigr) +
\ldots + \lambda _2\Bigl( p_{L (\varepsilon )} \bigl( S_{k} (\mathcal A ) \bigr) \Bigr)\leqslant
$$
$$
\leqslant
\alpha \bigl(S_1(\mathcal A)\bigr) + \ldots + \alpha \bigl(S_{k}(\mathcal A)\bigr) =
 \bigl( s_1^2 + \ldots + s_k^2 \bigr) \cdot \alpha  (\mathcal A).
$$
This contradicts the choice of $\varepsilon $.
Thus 
$\lambda _2 (p_\Pi \mathcal A) = 0$,
and the proof is finished.
\hfill $\square $

{\bf Proof of Theorem 2.}
Evidently, $(2)$ follows from $(3)$.
Also, $(2) $ easily implies $(1) $, see
\cite[p. 604, (1)]{Rushing-wild} or \cite[p. 136, (2)]{RS}.

$(1) \Rightarrow (2) $ is intuitively clear, but its rigorous proof
is rather long. 

\begin{figure}[!tbp]
  \begin{subfigure}[b]{0.31\textwidth}
  \centering
  \scalebox{0.65}{\includegraphics{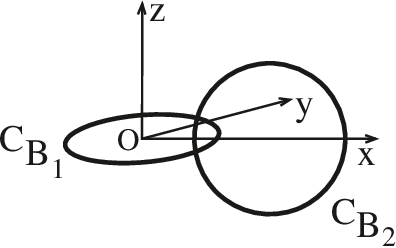}}
    \caption{Initial pair of circles.}
    \label{fig:fig3a}
  \end{subfigure}
  \begin{subfigure}[b]{0.31\textwidth}
   \centering
  \scalebox{0.65}{\includegraphics{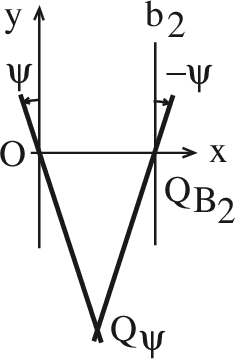}}
    \caption{Rotation scheme.}
    \label{fig:fig3b}
  \end{subfigure}
  \begin{subfigure}[b]{0.31\textwidth}
   \centering
  \scalebox{0.65}{\includegraphics{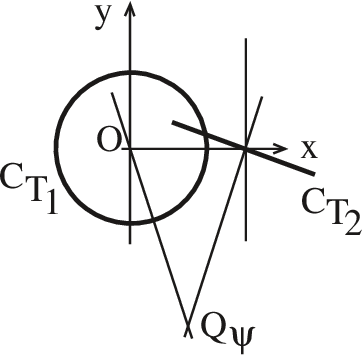}}
    \caption{After rotations.}
    \label{fig:fig3c}
  \end{subfigure}
  \caption{Illustration for the proof of the implication $(1)\Rightarrow (2)$ of Theorem 2.}
\end{figure}

\emph{Step 1: construct a special link.}
Take any number
$A \in \bigl( R_B+r_B , 2(R_B-r_B) \bigr)$;
the interval is non-empty since $R_B>3r_B$.
Define the circles $C_{B_1}, C_{B_2} \subset \mathbb R^3$
by their parametric equations:
$$
C_{B_1}: \quad
x=R_B \cos u, \ 
y=R_B \sin u, \ 
z=0
$$
and
$$
C_{B_2}: \quad
x=A+R_B \cos v, \ 
y=0, \ 
z=R_B \sin v ,
$$
where $u$ and $v$ run over the reals (see Fig. \ref{fig:fig3a}).

Let us show that the distance $d( C_{B_1}, C_{B_2}) > 2r_B$.
Take any points $P_1\in C_{B_1}$, $P_2\in C_{B_2}$; let 
$P_1 = (R_B \cos u , R_B \sin u , 0)$ and
$P_2 = (A+R_B \cos v , 0, R_B \sin v)$.
We have
$$
d^2 (P_1 , P_2) = (A+R_B \cos v - R_B\cos u)^2 +
(R_B \sin u)^2 + (R_B\sin v)^2 =
$$
$$
=
A^2 + 2AR_B \cos v - 2AR_B \cos u - 2R_B^2 \cos u \cos v + 2R_B^2 =: F(u,v) .
$$
Now
\[
 \left\{ \begin{array}{ll}
F'_u = 0 \\
F'_v = 0 \\
\end{array}
\right. \\
\iff
 \left\{ \begin{array}{ll}
(A+R_B\cos v) \sin u  = 0 \\
 (A-R_B\cos u) \sin v = 0 \\
\end{array}
\right. \\
\iff
 \left\{ \begin{array}{ll}
\sin u = 0 \\
\sin v = 0 \\
\end{array}
\right. \\
\]
(the last equivalence follows from the assumption $A>R_B+r_B > R_B$).
It remains to compute the distances $d(P_1,P_2)$ 
for $4$ pairs of points corresponding to $u,v \in \{ 0,\pi \}$;
the required inequality is easily obtained.

It is also easy to verify that $C_{B_1}$, $C_{B_2}$ are linked.

Now let $B_1, B_2\subset \mathbb R^3$ be closed $r_B$-neighborhoods
of $C_{B_1}$, $C_{B_2}$; they are
disjoint linked standard solid tori congruent to $B$.

\emph{Step 2: bend the link, complete the new link to a regular simple chain.} 
For
a straight line $\ell \| Oz$ and an angle $\psi $
we denote by $\rho_{\ell , \psi }$
the rotation transformation of
$\mathbb R^3$
with axis $\ell $ on angle $\psi $;
we rotate
in a counterclockwise direction when looking towards the plane $Oxy$ from the half-line  $Oz_{+}$ (the coordinate system $Oxyz$ is assumed to be a right-handed one).

Let $a_2 \| Oz$ and $b_2 \| Oy$ be the straight lines passing through $Q_{B_2}$.
For $\psi \in (0,\pi )$
let  $Q_{\psi } $ be the (unique)
common 
point of the lines 
$\rho _{Oz, \psi } (Oy) $ and $ \rho _{a_2, -\psi } (b_2)$, 
see Fig. \ref{fig:fig3b}.

By continuity, there is a number $\psi _0 \in (0,\pi )$
such that for each $\psi \in (0,\psi _0 )$ we have
$OQ_{\psi } > R_B + r_B$, and the set
$B_1\cap \rho _{a_2, -\psi } (B_2) = 
\rho _{Oz, \psi } (B_1) \cap \rho _{a_2, -\psi } (B_2) $ is empty.

Now let $m$ be any integer such that $m> \frac{\pi }{2 \psi _0}   $.
Fix the angle $\psi := \frac{\pi }{2m} < \psi _0$.
Solid tori 
$T_1:= B_1= \rho _{Oz, \psi } (B_1)$ and $T_2:=\rho _{a_2, -\psi } (B_2) $ 
are congruent to $B$.
(These tori depend on $\psi $, but we do not reflect this in the notation for brevity.)
By the above conditions, $T_1$ and $T_2$ are disjoint and linked.
Their central circles $C_{T_1} =  \rho _{Oz, \psi } (C_{B_1})$ and 
$C_{T_2}= \rho _{a_2, -\psi } (C_{B_2}) $, the way we see them from the axis 
$Oz_{+}$, are shown in Fig. \ref{fig:fig3c}.
For each $k=1,\ldots , m-1$
define
$$
T_{2k+1} := \rho _{\ell _\psi , -\frac{2\pi k}{m}} (T_1)
\quad
\text{ and }
\quad
T_{2k+2} := \rho _{\ell _\psi , -\frac{2\pi k}{m}} (T_2) ,
$$
where $\ell _\psi $ is the straight line passing through
$Q_\psi $ parallel to $Oz$.
Let $C_\psi $ be the circle in the $Oxy$-plane
having center $Q_\psi $ and radius $OQ_\psi $;
it contains the centers of $T_1,\ldots , T_{2m}$.
Define 
a standard solid torus $T$ as a closed 
$R_T$-neighborhood of $C_\psi $ for an arbirtary
$R_T \in \bigl( 0, OQ_\psi - (R_B+r_B)\bigr)$.
Now $T_1,\ldots , T_{2m}$
is a regular simple chain in $T$, as desired.

\begin{figure}[ht]
\centering
\scalebox{0.8}
{\includegraphics{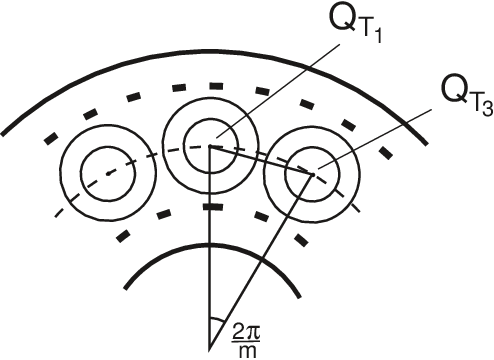}}
\caption{Illustration for the proof of the implication $(2)\Rightarrow (3)$ of Theorem 2.}
\label{fig:fig4}
\end{figure}

$(2) \Rightarrow (3) $. 
Let $T_1,\ldots , T_{2m}$ be a regular simple chain in $T$, where each $T_i$ is congruent to $B$;
let $\mathbf T$ 
be the standard solid torus similar to $B$ with
$C_{\mathbf T} = C_T$, see Fig. \ref{fig:fig4}; here two bold solid (resp. dashed) arcs denote the boundary of the solid torus $T$ (resp. $\mathbf T$), and 
the thin dashed arc represents the circle $C_{\mathbf T} = C_T$.
Assuming that $\sin \frac{\pi }{m} < \frac{r_B}{R_B}$,
we will prove that 
$r_{\mathbf T} > R_B+r_B $
(hence $T_i \subset \mathbf T$
for each $i=1,\ldots , 2m$).

Since $T_1$ and $T_3$ are disjoint,
the distance
$Q_{T_1}Q_{T_3} = 2R_T \sin \frac{\pi }{m} $
is greater than $\diam T_1 + \diam T_3 = 2(R_B+r_B)$. This together with the assumption $\sin \frac{\pi }{m} < \frac{r_B}{R_B}$
implies that
$$
R_T > \frac{R_B+r_B}{ \sin \frac{\pi }{m} } > 
\frac{R_B (R_B+r_B )}{ r_B} .
$$
From the 
similarity of $\mathbf T$ and $B$,
the equality $C_{\mathbf T} = C_T$, and the above inequality, we get
$$
r_{\mathbf T} = \frac{r_B R_{\mathbf T }}{R_B} = 
\frac{r_B R_{ T }}{R_B} > R_B+r_B .
$$
\hfill $\square $

\begin{figure}[ht]
\centering
\scalebox{0.6}
{\includegraphics{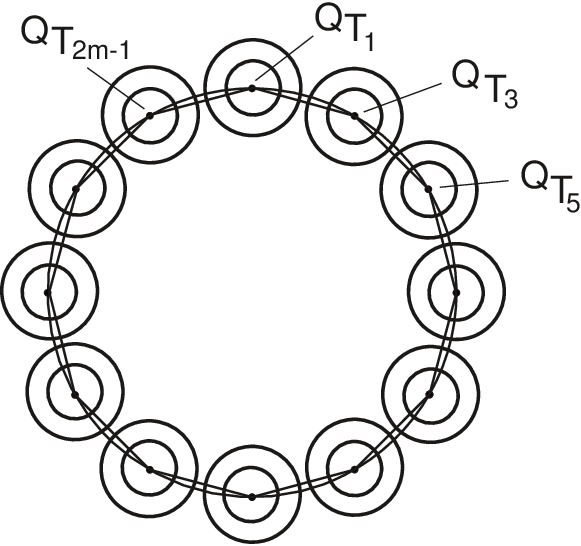}}
\caption{Illustration for the proof of part 1 of Theorem 3.}
\label{fig:fig5}
\end{figure}

{\bf Proof of Theorem 3.}
1) Intersect $T_1,T_3,\ldots , T_{2m-1}$ by the plane $\aff C_T$. We obtain circular rings 
with centers $Q_{T_1} , Q_{T_3}, \ldots , Q_{2m-1}\in C_T $;
the  radii of their outer circles equal $s(R_T+r_T)$, see  Fig. \ref{fig:fig5}.

Clearly $Q_{T_1}Q_{T_3}>2 s(R_T+r_T)$.
On the other side, 
the length of the closed broken line $Q_{T_1} Q_{T_3} \ldots  Q_{2m-1}$
is smaller than the length of the circle $C_T$, and
we get
$$m\cdot 2s (R_T+r_T) < m\cdot   Q_{T_1}Q_{T_3} 
=Q_{T_1} Q_{T_3} + Q_{T_3}Q_{T_5}+\ldots + 
 Q_{2m-1}Q_1
<  2\pi R_T .$$ 
Hence
$
ms<\frac{\pi R_T}{R_T+r_T} < \pi 
$
and
$
2ms^2 = 2s \cdot ms < 2\pi s $,
which proves 1).

2) The inclusion $T_1\subset T$ implies that
$s(R_T+r_T) \leqslant r_T$, and
$s\leqslant \frac{r_T}{R_T + r_T}$; this inequality was already noticed in \cite[p. 603]{Rushing-wild}, \cite[p. 136]{RS}, \cite[(2.1)]{Z}.
Now
$$
\frac{r_T}{R_T + r_T} < \frac{1}{2\pi } \iff 
\frac{r_T}{R_T} <  \frac{1}{2\pi - 1}  .
$$

3) We have 
$$ 2s(R_{T} + r_{T})  < Q_{T_1} Q_{T_3}  =
2R_T \sin \frac{\pi }{m} ,
$$
consequently
$$
s <  \frac{R_T \sin \frac{\pi }{m}}{ R_T +   r_T }
< \sin \frac{\pi }{m} <
 \frac{\pi }{m}.
$$
Since $2\pi ^2 <20$, we get
$\frac{\pi }{m} < \frac{1}{2\pi }$
for each  $m\geqslant 20$. 
\hfill $\square $

{\bf Proof of Theorem 4.}
Take any hyperplane $L\subset \mathbb R^n$.
The projection $p_L : \mathbb R^n\to L$
does not increase distances,
hence $m_{n-1} (p_L X) =0$.
By  the Szpilrajn theorem \cite[Thm. VII~3]{HW} we have $\dim p_L X \leqslant n-2$.

Now suppose that $\dim p_L X \leqslant n-3$; by
\cite{Wright}, \cite{Walsh-Wright} 
(for $n=3$ by \cite[{\bf 75}, p.~87--89]{Antoine-diss})
$X$ is tame, a contradiction.
This proves the theorem.
\hfill $\square $

\end{document}